\documentclass[12pt]{article}

\usepackage{hyperref}
\hypersetup{
  colorlinks   = true,    
  urlcolor     = blue,    
  linkcolor    = blue,    
  citecolor    = red      
}
\usepackage{epstopdf}
\usepackage[caption=false]{subfig}

\usepackage{amsmath}
\usepackage{csquotes}
\usepackage{tikz}
\usepackage{upgreek}
\usepackage{amssymb}
\usepackage{amsthm}
\usepackage{amsfonts}
\numberwithin{equation}{section}
\usepackage{booktabs,caption}
\usepackage[flushleft]{threeparttable}
\usepackage{stmaryrd}
\usepackage{wasysym}
\usepackage{makecell}
\usepackage{float}
\usepackage{graphicx}
\usepackage{lipsum}
\usepackage{orcidlink}
\usepackage{leftindex}

\newtheorem{theorem}{Theorem}[section]
\newtheorem{definition}[theorem]{Definition}

\newtheorem{lemma}[theorem]{Lemma}

\newtheorem{corollary}[theorem]{Corollary}
\newtheorem{notation}[theorem]{Notation}
\newtheorem{remark}[theorem]{Remark}

\theoremstyle{definition}
\newtheorem{example}[theorem]{Example}

\usepackage{nccmath}
\usepackage{mathtools}
\usepackage[reftex]{theoremref}
\usepackage{tikz-cd} 
\usepackage[left=1in, right=1in, top=1in, bottom=1in]{geometry}

\newtheorem{innercustomthm}{Theorem}
\newenvironment{mythm}[1]
  {\renewcommand\theinnercustomthm{#1}\innercustomthm}
  {\endinnercustomthm}

\begin{document}


\title{Inner and Outer Twisted Derivations of Cyclic Group Rings}

\author{Praveen Manju \orcidlink{0000-0001-6673-6123}
 and Rajendra Kumar Sharma \orcidlink{0000-0001-5666-4103}}
\date{}
\maketitle

\begin{center}
\noindent{\small Department of Mathematics, \\Indian Institute of Technology Delhi, \\ Hauz Khas, New Delhi-110016, India$^{1}$}
\end{center}

\footnotetext[1]{{\em E-mail addresses:} \url{praveenmanjuiitd@gmail.com}(Praveen Manju: Corresponding Author), \url{rksharmaiitd@gmail.com}(Rajendra Kumar Sharma).}

\medskip

\begin{abstract}
In this article, we study twisted derivations of cyclic group rings. Let $R$ be a commutative ring with unity, $G$ be a finite cyclic group, and ($\sigma, \tau$) be a pair of $R$-algebra endomorphisms of the group algebra $RG$, which are $R$-linear extensions of the group endomorphisms of $G$. In this article, we give two characterizations concerning $(\sigma, \tau)$-derivations of the group ring $RG$. First, we develop a necessary and sufficient condition for a $(\sigma, \tau)$-derivation of $RG$ to be inner. Second, we provide a necessary and sufficient condition for an $R$-linear map $D: RG \rightarrow RG$ with $D(1) = 0$ to be a $(\sigma, \tau)$-derivation. We also illustrate our theorems with the help of examples. As a consequence of these two characterizations, we answer the well-known twisted derivation problem for $RG$: Under what conditions are all $(\sigma, \tau)$-derivations of $RG$ inner? Or is the space of outer $(\sigma, \tau)$-derivations trivial? More precisely, we give a sufficient condition under which all $(\sigma, \tau)$-derivations of $RG$ are inner and a sufficient condition under which $RG$ has non-trivial outer $(\sigma, \tau)$-derivations. Our result helps in generating several examples of non-trivial outer derivations.
\end{abstract}



\section{Introduction}\label{section 1}
Derivations play an essential role in mathematics and physics. The theory of derivations has been developed in rings and numerous algebras and helps to study and understand their structure. For example, von Neumann algebras \cite{Brear1992}, MV-algebras \cite{KamaliArdakani2013, Mustafa2013}, BCI-algebras \cite{Muhiuddin2012}, Banach algebras \cite{Raza2016}, incline algebras that have many applications \cite{Kim2014}, lattices that are very important in fields such as information theory: information recovery, information access management, and cryptanalysis \cite{Chaudhry2011}. Operator algebras, $\mathbb{C}^{*}$-algebras, differentiable manifolds, and representation theory of Lie groups are being studied using derivations \cite{Klimek2021}. For a historical account and further applications of derivations, we refer the reader to \cite{Atteya2019, Haetinger2011, MohammadAshraf2006}. Jacobson \cite{Jacobson1964} introduced the idea of an $(s_{1}, s_{2})$-derivation, which was later commonly called the $(\sigma, \tau)$ or $(\theta, \phi)$-derivation. These twisted derivations have been studied extensively in prime and semiprime rings and have been mainly used in solving functional equations \cite{Brear1992}. Twisted derivations have huge applications. They are used in multiplicative deformations and discretizations of derivatives that have many applications in models of quantum phenomena and the analysis of complex systems and processes. They are extensively investigated in physics and engineering. Using twisted derivations, Lie algebras are generalized to hom-Lie algebras, and the central extension theory is developed for hom-Lie algebras analogously to that for Lie algebras. Just as Lie algebras were initially studied as algebras of derivations, Hom-Lie algebras were defined as algebras of twisted derivations. The generalizations (deformations and analogs) of the Witt algebra, the complex Lie algebra of derivations on the algebra of Laurent polynomials $\mathbb{C}[t, t^{-1}]$ in one variable, are obtained using twisted derivations. Deformed Witt and Virasoro-type algebras have applications in analysis, numerical mathematics, algebraic geometry, arithmetic geometry, number theory, and physics. We refer the reader to \cite{Hartwig2006, Ilwale2023, Larsson2017, Siebert1996} for details. Twisted derivations have been used to generalize Galois theory over division rings and in the study of $q$-difference operators in number theory. Derivations, especially $(\sigma, \tau)$-derivations of rings and group rings, have various applications in coding theory \cite{Creedon2019, Boucher2014}. For more applications of $(\sigma, \tau)$-derivations, we refer the reader to \cite{AleksandrAlekseev2020} and the references within.

Sergei D. Silvestrov, in his several papers  \cite{abramov20203, arfa2019classification, armakan2020non, back2018hom, elchinger2016brackets, elhamdadi2022derivation, farhangdoost2023simply, Hartwig2006, laraiedh2021transposed, larsson2006some, larsson2008quasi, larssonq, larsson2005lie, larsson2005quasi, larsson2006quasi, larsson2007quasi, ma2017rota, ma2021rota, mabrouk2023pseudo, mabrouk2021generalized, musonda2020twisted, oinert2013maximal, richard2009note, richard2008quasi, unknown, saadaoui2019lambda} has done a lot of work on twisted derivations of algebraic structures. In \cite{abramov20203}, the authors constructed $3$-Hom-Lie algebras using an associative commutative algebra with a twisted derivation and involution. In \cite{arfa2019classification}, the authors studied the classification of multiplicative two-dimensional Hom-Leibniz algebras and investigated derivations of Hom-Leibniz algebras and superalgebras. In \cite{armakan2020non}, the authors investigated some important basic properties of simple Hom-Lie superalgebras. In doing this, the authors employed derivations of these Hom-Lie superalgebras.
In \cite{back2018hom}, the authors studied Hom-associative Ore extensions. Some existence theorems employed the use of twisted derivations. In \cite{elchinger2016brackets}, the authors studied some brackets defined on $(\tau, \sigma)$-derivations satisfying quasi-Lie identities. The authors defined quasi-Lie and Hom-Lie algebras in the case of $(\tau, \sigma)$-derivations. In \cite{elhamdadi2022derivation}, the authors studied the derivations of quandle algebras. which were introduced in an analogous way to the theory of group rings. In \cite{farhangdoost2023simply}, the authors considered complete Hom-Lie superalgebras. They obtained conditions under which the linear space of $\alpha^{s}$-derivations of a Hom-Lie superalgebra is complete and simply complete. In \cite{laraiedh2021transposed}, the authors introduced the notion of transposed Hom-Poisson and Hom-pre-Lie Poisson algebras using derivations (see, for example, Proposition 2.19, Proposition 2.23). In \cite{larsson2006some}, the authors obtained multiparameter families of almost quadratic algebras. They explored the general quasi-deformation scheme in the particular cases of $\mathbb{F}$-algebras ($\mathbb{F}$ a field of characteristic $0$) as applied to $\mathfrak{sl}_{2}(\mathbb{F})$ in the particular cases of $\mathcal{A} = \mathbb{F}[t]$ and $\mathcal{A} = \frac{\mathbb{F}[t]}{\langle t^{4} \rangle}$. Note that deformation theory is defined using twisted derivations. In \cite{larsson2008quasi}, the authors showed that one could obtain interesting new multi-parameter families of almost quadratic algebras when the so-called quasi-deformation method is applied to $\mathfrak{sl}_{2}(\mathbb{F})$ via representations by twisted derivations on the algebra $\frac{\mathbb{F}[t]}{\langle t^{N} \rangle}$, where $\mathbb{F}$ is a field of characteristic $0$. In \cite{larssonq}, the authors generalized a result of V. Abramov on $q$-differential graded algebras and showed its relation to $N$-complexes in explicit terms. In doing that, the authors made a significant use of twisted derivations. In \cite{larsson2005lie}, the authors studied the deformation of $\mathfrak{sl}_{2}(\mathbb{F})$ and discussed some particular algebras appearing in this deformation scheme. Here $\mathbb{F}$ is a field of characteristic $0$. In \cite{larsson2005quasi}, the authors introduced the notion of quasi-Hom-Lie algebras, a new class of deformations of Lie algebras associated with twisted derivations. Quasi-Hom-Lie-algebras include Hom-Lie algebras, color Lie algebras, Lie superalgebras and Lie algebras, as well as other more exotic types of algebras, which then can be viewed as a kind of deformation of Lie algebras in some larger category. In \cite{larsson2006quasi}, the authors constructed quasi-deformation of $\mathfrak{sl}_{2}(\mathbb{R})$ using more general representations of $\mathfrak{sl}_{2}(\mathbb{R})$ in terms of derivations on a suitable function algebra $\mathcal{A}$ in an indeterminate $t$, and then specializing to the case $\mathbb{R}[t,t^{-1}]$ where derivations are replaced by twisted derivations. In \cite{larsson2007quasi}, the authors obtained quasi-deformations of the simple $3$-dimensional Lie algebra $\mathfrak{sl}_{2}(\mathbb{F})$ ($\mathbb{F}$ a field of characteristic $0$) by using a deformation method that relies heavily on twisted derivations. In \cite{ma2017rota}, the authors studied covariant bialgebras via coderivations. They gave an alternative definition of covariant bialgebras through two coderivations and provided a characterization. In \cite{ma2021rota}, the authors introduced a new type of bialgebra, namely, mixed bialgebra (as a generalization of infinitesimal bialgebra) which consists of an associative algebra and a coassociative coalgebra satisfying the compatible condition determined by two coderivations. In \cite{mabrouk2023pseudo}, the authors gave the concept of pseudo-Euclidean Hom-alternative superalgebras and symplectic Hom-alternative superalgebra and the construction of symplectic Hom-alternative superalgebras from pseudo-Euclidean symplectic Hom-alternative superalgebras and an antisymmetric superderivation. In \cite{mabrouk2021generalized}, the authors studied the notion of generalized derivations and Rota-Baxter operators of $n$-ary Hom-Nambu and $n$-Hom-Lie superalgebras and their relation with generalized derivations and Rota-Baxter operators of Hom-Lie superalgebras. In \cite{musonda2020twisted}, it was shown that there exist twisted derivations or twisted derivation-like operators generalizing the Jackson derivative that, together with the multiplication operator, gives a concrete representation of a certain multi-parameter family of unital associative complex algebras generated by certain elements satisfying some commutation relations. In \cite{oinert2013maximal}, the authors studied the simplicity of Ore extension rings $R[x; \sigma, \delta]$ ($x \notin R$ such that $1, x, x^{2}, ...$ form a basis of $R[x; \sigma, \delta]$ as a left $R$-module, $\sigma:R \rightarrow R$ ring endomorphism, $\delta:R \rightarrow R$ a $\sigma$-derivation and for all $r \in R$, $xr = \sigma(r)x + \delta(r)$). Many results were obtained in terms of the twisted derivation $\delta$ of $R$ (see, for example, Proposition 3.11, Remark 4.1, Theorem 4.15). In \cite{richard2009note}, the authors studied the properties of the bracket defined by $\sigma$-derivations of $\mathbb{C}[t_{1}^{\pm 1}, ..., t_{n}^{\pm 1}]$ and emphasized on the question to know when this bracket defines a Hom-Lie structure. In \cite{richard2008quasi}, the authors proved that a certain bracket on the $\sigma$-derivations of a commutative algebra preserves inner derivations (see Proposition 2.3.1), and based on this obtained structural results providing new insights into $\sigma$-derivations on Laurent polynomials in one variable, $\mathbb{C}[t, t^{-1}]$ which is a UFD. In \cite{Hartwig2006}, the authors developed an approach to deformations of the Witt and Virasoro algebras based on twisted derivations. The authors proved that if $\mathcal{A}$ (an algebra over $\mathbb{C}$) is a unique factorization domain (UFD), then the $\mathcal{A}$-module $\mathcal{D}_{(\sigma, id_{\mathcal{A}})}(\mathcal{A})$ of all $(\sigma, id_{\mathcal{A}})$-derivations of $\mathcal{A}$ is a free $\mathcal{A}$-module of rank one. Here $\sigma$ and $\tau$ are two different algebra endomorphisms of $\mathcal{A}$ and the generator is a $(\sigma, \tau)$-derivation given by $\Delta = \frac{\tau - \sigma}{g}$, where $g = \text{gcd}((\tau - \sigma)(\mathcal{A}))$ (see Theorem 4). The authors also gave the notion of generalized Witt algebra using a twisted derivation. In \cite{richard2008quasi}, the authors further noted that a $(\sigma, id_{\mathcal{A}})$-derivation $a \Delta$ is inner if and only if $g$ divides $a$. More precisely, the $\mathcal{A}$-submodule of all inner $(\sigma, id_{\mathcal{A}})$-derivations is of rank one with generator $g \mathcal{A} \Delta$ (see Proposition 2.4.1). Then the authors obtained some more deep and precise results for the algebra $\mathcal{A} = \mathbb{C}[t, t^{-1}]$, which is also a UFD. Note that, in this article, we consider our set of $(\sigma, \tau)$-derivations of commutative group ring $RG$ ($R$ a commutative unital ring, $G$ a finite cyclic group) as an $R$-module rather than an $RG$-module. In \cite{unknown}, the authors considered the Ore extension algebra for the algebra $\mathcal{A}$ of functions with finite support on a countable set. They derived explicit formulas for twisted derivations on $\mathcal{A}$. They gave a description of twisted derivations on the algebra of functions on a finite set from which it is observed that there are no non-trivial derivations on $\mathbb{R}^{n}$. In \cite{saadaoui2019lambda}, the authors generalized the results about generalized derivations of Lie algebras to the case of BiHom-Lie algebras. In \cite{silvestrov2007paradigm}, the author reviewed the main constructions and examples of quasi-deformations of Lie algebras via twisted derivations leading to quasi-Lie algebras. In \cite{larsson2014arithmetic}, the author studied Hom-Lie algebras and global-equivariant Hom-Lie algebras in the context of twisted derivations. The authors obtained some nice results on twisted module derivations (see, for example, Section 2.2: Lemma 2.4, Proposition 2.5, Theorem 2.5) and global twisted derivations (see Section 2.3).

This article considers the $(\sigma, \tau)$-derivations of the pure algebraic structure, namely, the group ring. For a history of group rings, we refer the reader to {\cite[Chapter 3]{CPM2002}}.
Derivations and $(\sigma, \tau)$-derivations of group rings (defined purely algebraically) have not received much attention. The study of derivations of group rings begins with the paper \cite{Smith1978}, where the author studies derivations of group rings of a finitely-generated, torsion-free, nilpotent group $G$ over a field $\mathbb{F}$. For instance, it is shown that such group rings always contain an outer derivation. She, using the notion of derivations of $\mathbb{F}G$, proved that the Krull dimension of such a group ring $\mathbb{F}G$ is finite and equals $\mathbb{F}$-rank of $G$ if $G$ satisfies certain additional conditions. In {\cite[Theorem 1]{Spiegel1994}}, the main theorem of \cite{Spiegel1994}, it is proved that every derivation of an integral group ring of a finite group is inner. The author also described the derivations of $\mathbb{Z}G$ for any abelian group $G$. In {\cite[Theorem 1.1]{Ferrero1995}}, the main theorem of \cite{Ferrero1995}, the authors demonstrate that if $G$ is a torsion group whose center $Z(G)$ has finite index in $G$ and $R$ is a semiprime ring such that $\text{char}(R)$ is either $0$ or does not divide the order of any element of $G$, then every $R$-derivation of the group ring $RG$ is inner. In \cite{Chaudhuri2019} and \cite{Chaudhuri2021}, the author generalizes the above results of \cite{Spiegel1994} and \cite{Ferrero1995} respectively to $(\sigma, \tau)$-derivations of group rings of finite groups over a field and an integral domain where $\sigma, \tau$ satisfy certain conditions. In \cite{Chaudhuri2019}, the author proves the following main theorem:
\begin{theorem}[{\cite[Theorem 1.1]{Chaudhuri2019}}]
Let $G$ be a finite group and $R$ be an integral domain with unity such that $|G|$ is invertible in $R$. Let $\sigma, \tau$ be $R$-algebra endomorphisms of $RG$ such that they fix $Z(RG)$ elementwise.
\begin{enumerate}
\item[(i)] If $R$ is a field, then every $(\sigma, \tau)$-derivation of $RG$ is inner.
\item[(ii)] If $R$ is an integral domain which is not a field and if $\sigma, \tau$ are $R$-linear extensions of the group homomorphisms of $G$, then every $(\sigma, \tau)$-derivation of $RG$ is inner.
\end{enumerate}
\end{theorem} 
\noindent She also gives the following as a corollary of the above theorem. If $G$ is a finite group, $R = \mathbb{Z}$ and $\sigma, \tau$ satisfy the conditions of part (ii) of the above theorem, then every $(\sigma, \tau)$-derivation of $\mathbb{Z}G$ is inner. Below, we compare the work done in \cite{Chaudhuri2019} in the form of the above theorem with the work we have done in this article to ensure the novelty of our article. Also, our techniques are different. \begin{enumerate}
\item[•] In the above theorem, the author assumes that the coefficient ring $R$ is either a field or an integral domain, but in our work of this article, we take $R$ to be an arbitrary commutative ring with unity. 
\item[•] Further, in the above theorem, the authors takes $\sigma$ and $\tau$ to be the $R$-algebra endomorphisms of $RG$ which must fix the center $Z(RG)$ of the group ring $RG$ elementwise, and if $R$ is an integral domain $\sigma, \tau$ are also $R$-linear extensions of the group endomorphisms of $G$. But, in our work of this article, we do not make the above assumption that $\sigma, \tau$ fix $Z(RG)$ elementwise. In fact, we take $\sigma$ and $\tau$ to be $R$-algebra endomorphisms of the cyclic group ring $RG$ of a cyclic group $G$ over a commutative unital ring $R$, which are $R$-linear extensions of the group endomorphisms of $G$; note that such a $\sigma, \tau$ may not fix $Z(RG)$ elementwise. In fact, only the identity group endomorphism fixes $Z(RG)$ elementwise.
\item[•] Most importantly, the above theorem does not give sufficient conditions for group rings to have outer derivations. But, in our article, we are studying the twisted derivation problem in which we have determined the sufficient condition under which $RG$ has non-trivial outer derivations and the condition under which all derivations of $RG$ are inner.
\end{enumerate} 
In {\cite[Theorem 1.1]{Chaudhuri2021}}, the main theorem of \cite{Chaudhuri2021}, the author proved that if $R$ is a semiprime ring with unity such that either $R$ does not have torsion elements or that if $R$ has $p$-torsion elements, then $p$ does not divide $|G|$, where $G$ is a torsion group such that $[G:Z(G)] < \infty$ and if $\sigma, \tau$ are $R$-algebra endomorphisms of $RG$ which fix $Z(RG)$ elementwise, then there is a ring $T$ such that $R \subseteq T$, $Z(R) \subseteq Z(T)$ and for the natural extensions of $\sigma, \tau$ to $TG$, $H^{1}(TG, \leftindex_{\sigma} {TG}_{\tau}) = \{0\}$, the degree $1$ Hochschild cohomology, that is, there exists a ring extension $T$ of the semiprime ring $R$ such that all $(\sigma, \tau)$-derivations of $TG$ are inner for the natural extensions $\sigma, \tau$ to $TG$. Also, in {\cite[Theorem 4.6]{Chaudhuri2021}}, the author gives an element in $RG$ (a commutative group algebra of a commutative unital ring $R$ with certain conditions over an abelian group $G$) that makes a $(\sigma, \tau)$-derivation of $RG$ inner under the condition that there exists some $b \in RG$ with $(\tau - \sigma)(b)$ invertible in $RG$, where $\sigma, \tau$ are $R$-algebra endomorphisms of $RG$.

A. A. Arutyunov, in his several papers (\cite{AleksandrAlekseev2020, Arutyunov2021, A.A.Arutyunov2020, Arutyunov2023, Arutyunov2020a, Arutyunov2020, arutyunov2019smooth}), studies derivations using topology and characters. In \cite{AleksandrAlekseev2020}, the authors considered the $(\sigma, \tau)$ and $(\sigma, Id)$-derivations of the group algebra $\mathbb{C}G$ of a discrete countable group $G$, where $\sigma, \tau$ are $\mathbb{C}$-algebra endomorphisms of $\mathbb{C}G$ which are $\mathbb{C}$-linear extensions of the group endomorphisms of $G$. In {\cite[Theorem 4]{AleksandrAlekseev2020}}, the authors proved a decomposition theorem expressing the $\mathbb{C}$-vector space of all $(\sigma, \tau)$-derivations of $\mathbb{C}G$ of a finitely generated $(\sigma, \tau)$-FC group, where a $(\sigma, \tau)$-FC group is a group in which every $(\sigma, \tau)$-conjugacy class is of finite size. As a corollary, in {\cite[Corollary 5]{AleksandrAlekseev2020}}, it was concluded that every $(\sigma, \tau)$-derivation of the group algebra $\mathbb{C}G$ of a finite group $G$ is inner, where $\sigma, \tau$ are as above, $\mathbb{C}$-algebra endomorphisms of $\mathbb{C}G$ which are $\mathbb{C}$-linear extensions of the group endomorphisms of $G$. But, in this article, we work over an arbitrary commutative unital ring $R$ instead of a field. Also, our techniques are different.

In \cite{Arutyunov2021}, the authors studied the ordinary derivations of finite and FC groups using topological and character techniques. They obtained decomposition theorems for derivations of the group ring $AG$ of a finite group or FC group (a group in which every conjugacy class is of finite size) $G$ over a commutative unital ring $A$ {\cite[Theorems 3.2, 4.1, 4.3; Corollary 4.2.4]{Arutyunov2021}}. More precisely, the authors proved that for a finitely generated FC-group $G$, the $A$-module of all derivations of $AG$ is isomorphic to the direct sum of its $A$-submodule of all inner derivations of $AG$ and $\bigoplus_{[u] \in G^{G}} \text{Hom}_{Ab}(Z(u),A),$ where $\text{Hom}_{Ab}(Z(u),A)$ is the set of additive homomorphisms from the centralizer $Z(u)$ of a fixed element $u \in G$ to the ring $A$ and $G^{G}$ is the set of all possible distinct conjugacy classes $[u]$ ($u \in G$) in $G$ {\cite[Theorems 3.2]{Arutyunov2021}}. Furthermore, if $G$ is a finite group, then the $A$-module of all outer derivations of $AG$ is isomorphic to $\bigoplus_{[u] \in G^{G}} \text{Hom}_{Ab}(Z(u),A)$ {\cite[Theorems 4.1]{Arutyunov2021}}. In {\cite[Theorems 4.3]{Arutyunov2021}}, the authors gave a characterization under which all ordinary derivations of $AG$ are inner for a finite group $G$ and a finite ring $A$. In {\cite[Corollary 4.2.4]{Arutyunov2021}}, the authors proved that for a torsion-free ring $A$ and a finite group $G$, all derivations are inner. In our article, we are considering more generalized $(\sigma, \tau)$-derivations of a cyclic group ring $RG$ over a commutative unital ring $R$ with no further conditions on it like finiteness and torsion-free.

In \cite{A.A.Arutyunov2020} also, a description of derivations is given in terms of characters of the groupoid of the adjoint action of the group. A method of the description of the space of all outer derivations of the complex group algebra $\mathbb{C}G$ of a finitely presentable discrete group $G$ is given. In \cite{Arutyunov2023}, the author studied derivations of a complex group algebra $\mathbb{C}G$ of a finitely generated group $G$. The author constructed the ideals of inner and quasi-inner derivations and established a connection between derivations and characters on the groupoid of the adjoint action (see Proposition 2.11 and Corollary 3.2). The author gave a description of the space of outer and quasi-outer derivations of $\mathbb{C}G$ using the methods of combinatorial group theory, in particular, the number of ends of the group $G$ and the number of ends of the connected components of a conjugacy diagram (see, for example, Theorem 3.6, Corollary 3.7, Theorem 4.2, Corollary 4.3, Corollary 4.5, Corollary 4.7, Proposition 5.1). In \cite{Arutyunov2020a}, the authors gave a description of the space of all inner and outer derivations of a group algebra of a finitely presented discrete group in terms of the character spaces of the 2-groupoid of the adjoint action of the group. In \cite{Arutyunov2020b}, the author established a connection between derivations of complex group algebras and the theory that studies the ends of topological spaces. In {\cite[Theorems 5 and 6]{Arutyunov2020b}}, the author gave a description of the space of outer and quasi-outer derivations of a complex group algebra.

In \cite{Arutyunov2020}, the author studied derivation spaces in the group algebra $\mathbb{C}G$ of a generally infinite non-commutative discrete group $G$ in terms of characters on a groupoid associated with the group. In {\cite[Theorem 2]{Arutyunov2020}}, the author constructed a subalgebra of non-inner derivations of $\mathbb{C}G$, which can be embedded in the algebra of all outer derivations. In {\cite[Theorem 3]{Arutyunov2020}}, the author obtained the necessary conditions under which a character defines a derivation. In {\cite[Proposition 8]{Arutyunov2020}}, the author described the space of all derivations of a complex group algebra of a free abelian group of finite rank. In {\cite[Theorem 4 and Theorem 5]{Arutyunov2020}}, the author gave an explicit description of the space of all derivations and outer derivations of the complex group algebra of a rank 2 nilpotent group.
In \cite{arutyunov2019smooth}, the authors gave the description of the space of all outer derivations of the group algebra $\mathbb{C}G$ of a finitely presented discrete group in terms of the Cayley complex of the groupoid of the adjoint action of the group (see Theorem 3, Theorem 4, Corollary 1, Corollary 2). In {\cite[Section 5]{arutyunov2019smooth}}, the authors described the space of ordinary derivations of a complex group algebra of the additive group $\mathbb{Z}$, a finitely generated free abelian group and a finitely generated free group. They showed that the space of outer derivations is isomorphic to the one-dimensional compactly supported group of the Cayley complex over $\mathbb{C}$. In \cite{mishchenko2020description}, the author studied the derivations of a complex group algebra $\mathbb{C}G$ of a finitely generated discrete group. In {\cite[Section 3]{mishchenko2020description}}, the author gave the description of the derivations of $\mathbb{C}G$ using characters on the groupoid of the adjoint action of the group (see Theorems 2 and 3). In {\cite[Section 4]{mishchenko2020description}}, the author gave a description of the space of outer derivations by establishing the isomorphism of this space with the one-dimensional cohomology of the Cayley complex of the groupoid (see Corollary 1). In {\cite[Section 5]{mishchenko2020description}}, the author described the space of ordinary derivations of a complex group algebra of the additive group $\mathbb{Z}$, a finitely generated free abelian group and a finitely generated free group.

In \cite{OrestD.Artemovych2020}, the authors considered the group ring of a group $G$ over a unital ring $R$ such that all prime divisors of orders of elements in $G$ are invertible in $R$. They proved that if $R$ is finite and $G$ is a torsion FC-group, then all derivations of $RG$ are inner. They obtained similar results for other classes of groups $G$ and rings $R$. In \cite{Creedon2019}, the authors studied the derivations of group rings over a commutative unital ring $R$ in terms of the generators and relators of the group. In {\cite[Theorem 2.5]{Creedon2019}}, which is the main theorem of their article, the authors gave a necessary and sufficient condition under which a map from the generating set of the group $G$ to the group ring $RG$ can be extended to a derivation of $RG$. 
Applying this characterization, the authors in {\cite[Theorem 3.4]{Creedon2019}} classified the ordinary derivations of commutative group algebras over a field of prime characteristic $p$ by giving the dimension and a basis for the vector space of all derivations. In \cite{Manju2023a}, the authors classified all $(\sigma, \sigma)$-derivations of a commutative group algebra $\mathbb{F}G$ of an abelian group $G$ over a field $\mathbb{F}$ of positive characteristic, where $\sigma$ is an $\mathbb{F}$-algebra endomorphism of $\mathbb{F}G$ which is an $\mathbb{F}$-linear extension of a group endomorphism of $G$. But, in this article, we consider the general case of $(\sigma, \tau)$ rather than $(\sigma, \sigma)$ for finite cyclic groups, and also, we consider our coefficient ring to be commutative unital ring unlike above results which considered the coefficient ring to be a field.

This manuscript presents a new but simple approach to studying twisted derivations in group rings. Also, the ring $R$ is taken, in general, as a commutative unital ring with no other restriction, unlike most earlier articles, which have taken $R$ as a field, an integral domain, or a ring with other conditions on its characteristic. The manuscript has been divided into five sections. Section \ref{section 2} states some preliminaries containing some basic definitions and useful results that we will need later. This article considers the twisted derivation problem: Are all twisted derivations inner? Or is the space of outer twisted derivations trivial? We refer the reader to \cite{AleksandrAlekseev2020, A.A.Arutyunov2020, Arutyunov2020b, Arutyunov2020} for the history and importance of the derivation problem. This article considers the analogous problem for $(\sigma, \tau)$-derivations of the cyclic group rings. Let $R$ be a commutative unital ring, $G = \langle x \mid x^{n} = 1\rangle$ be a finite cyclic group of order $n$, and $(\sigma, \tau)$ be a pair of $R$-algebra endomorphisms of $RG$, which are $R$-linear extensions of the group endomorphisms of $G$. Section \ref{section 4} gives a necessary and sufficient condition under which a $(\sigma, \tau)$-derivation of $RG$ is inner. The main theorem of Section \ref{section 4} is stated as follows and has been proved with the help of four lemmas, making extensive and beautiful use of linear algebra over commutative unital rings.
\begin{mythm}{3.5}
Let $R$ be a commutative unital ring, $G = \langle x \mid x^{n} = 1\rangle$ be a cyclic group of order $n$, and $(\sigma, \tau)$ be a pair of $R$-algebra endomorphisms of $RG$ which are $R$-linear extensions of the group endomorphisms of $G$ (so that $\sigma(x) = x^{u}$ and $\tau(x) = x^{u+v}$ for some $u, v \in \{0, 1, ..., n-1\}$). Let $D:RG \rightarrow RG$ be a $(\sigma, \tau)$-derivation. Suppose that $D(x) = \sum_{i=0}^{n-1}c_{i}x^{i}$, $d = \text{gcd}(v,n)$, and $m = \frac{n}{d}$. Then $D$ is inner if and only if the following $d$ equations hold simultaneously: \begin{equation*}
\sum_{j=0}^{m - 1} c_{i+jd} = 0, ~~~~~ i \in \{0, 1, ..., d-1\}.
\end{equation*}
\end{mythm}

We first see the motivation behind Lemma \ref{lemma 3.1}. Given a commutative ring $R$ with unity, a group $G$, and a pair $(\sigma, \tau)$ of $R$-algebra endomorphisms of $RG$, a $(\sigma, \tau)$-derivation $D$ satisfies the identity \begin{equation*} D(\alpha^{k}) = \left(\sum_{i+j=k-1} \sigma(\alpha^{i}) \tau(\alpha^{j})\right) D(\alpha)\end{equation*} for all $\alpha \in RG$ and for all $k \in \mathbb{N}$, where $i,j$ run over non-negative integers. A natural question one may ask is: under what assumptions on an $R$-linear map $D: RG \rightarrow RG$ satisfying the above relations is a $(\sigma, \tau)$-derivation? Similar questions have been posed for derivations of a ring (for example, see \cite{Bridges1984, Vukman2005}). In this article, we consider the question for the $(\sigma, \tau)$-derivations of a cyclic group ring $RG$, that is, when $R$ is a commutative ring with unity, $G$ is a finite cyclic group $G = \langle x \mid x^{n}=1\rangle$ of order $n$ and $\sigma, \tau$ are $R$-algebra endomorphisms of $RG$, which are $R$-linear extensions of the group endomorphisms of $G$. Note that $\sigma(x) = x^{u}$ and $\tau(x) = x^{u+v}$ for some $u, v \in \{0, 1, ..., n-1\}$. In Lemma \ref{lemma 3.1}, we give a characterization under which an $R$-linear map $D:RG \rightarrow RG$ with $D(1)=0$ is a $(\sigma, \tau)$-derivation. We also see its useful consequences in the form of corollaries. Not only this, our \th\ref{lemma 3.1} is very crucial in answering the twisted derivation problem. As a consequence of \th\ref{lemma 3.1}, \th\ref{theorem 4.1}, and Remark \ref{remark 2.5}, we answer the twisted derivation problem for $RG$ in the most important result of this manuscript, namely, \th\ref{theorem 3.4}, stated below:
\begin{mythm}{4.4}
Let $R$ be a commutative unital ring, $G = \langle x \mid x^{n} = 1\rangle$ be a cyclic group of order $n$, and $(\sigma, \tau)$ be a pair of $R$-algebra endomorphisms of $RG$, which are $R$-linear extensions of the group endomorphisms of $G$ (so that $\sigma(x) = x^{u}$ and $\tau(x) = x^{u+v}$ for some $u, v \in \{0, 1, ..., n-1\}$). Then the following statements hold:
\begin{enumerate}
\item[(i)] If $\text{gcd}(v,n)$ is invertible in $R$, then every $(\sigma, \tau)$-derivation of $RG$ is inner, that is, $\mathcal{D}_{(\sigma, \tau)}(RG) = \text{Inn}_{(\sigma, \tau)}(RG)$.
\item[(ii)] If $\text{char}(R)$ divides $\text{gcd}(v,n)$, then $\text{Inn}_{(\sigma, \tau)}(RG) \subsetneq \mathcal{D}_{(\sigma, \tau)}(RG)$, that is, $RG$ has non-zero outer $(\sigma, \tau)$-derivations.
\end{enumerate}
\end{mythm}

In Section \ref{section 5}, we give some examples illustrating \th\ref{theorem 4.1} and \th\ref{theorem 3.4}.

\section{Preliminaries}\label{section 2}
Below, we state some basic definitions. Unless otherwise stated, $R$ denotes a commutative unital ring and $\mathcal{A}$ an associative $R$-algebra. Let $(\sigma, \tau)$ be a pair of $R$-algebra endomorphisms of $\mathcal{A}$. 

\begin{definition}\label{definition 2.1}
An $R$-linear map $d:\mathcal{A} \rightarrow \mathcal{A}$ that satisfies $d(\alpha \beta) = d(\alpha) \beta + \alpha d(\beta)$ for all $\alpha, \beta \in \mathcal{A}$, is called a derivation of $\mathcal{A}$. It is called inner if there exists some $\beta \in \mathcal{A}$ such that $d(\alpha) = \beta \alpha - \alpha \beta$ for all $\alpha \in \mathcal{A}$, and then we denote it by $d_{\beta}$. The elements of the quotient of the $R$-module of all derivations of $\mathcal{A}$ by the $R$-submodule of all inner derivations are called outer derivations.
\end{definition}

\begin{definition}\label{definition 2.2}
A $(\sigma, \tau)$-derivation $D:\mathcal{A} \rightarrow \mathcal{A}$ is an $R$-linear map that satisfies the $(\sigma, \tau)$-twisted generalized identity: $D(\alpha \beta) = D(\alpha)\tau(\beta) + \sigma(\alpha) D(\beta)$ for all $\alpha, \beta \in \mathcal{A}$. It is called inner if there exists some $\beta \in \mathcal{A}$ such that $D(\alpha) = \beta \tau(\alpha) - \sigma(\alpha) \beta$ for all $\alpha \in \mathcal{A}$, and then we denote it by $D_{\beta}$. The elements of the quotient of the $R$-module of all $(\sigma, \tau)$-derivations of $\mathcal{A}$ by the $R$-submodule of all inner $(\sigma, \tau)$-derivations are called outer $(\sigma, \tau)$-derivations.
\end{definition}

\begin{remark}\label{remark 2.3}
If $\sigma = \tau = id_{\mathcal{A}}$ (the identity map on $\mathcal{A}$), then the usual Leibniz identity holds, and $(\sigma, \tau)$-derivation and inner $(\sigma, \tau)$-derivation respectively become the ordinary derivation and ordinary inner derivation of $\mathcal{A}$. We denote the set of all $(\sigma, \tau)$-derivations of $\mathcal{A}$ by $\mathcal{D}_{(\sigma, \tau)}(\mathcal{A})$, the set of all inner $(\sigma, \tau)$-derivations of $\mathcal{A}$ by $\text{Inn}_{(\sigma, \tau)}(\mathcal{A})$ and the set of all outer $(\sigma, \tau)$-derivations of $\mathcal{A}$ by $\text{Out}_{(\sigma, \tau)}(\mathcal{A})$. Defining componentwise sum and module action, $\mathcal{D}_{(\sigma, \tau)}(\mathcal{A})$ becomes an $R$- as well as $\mathcal{A}$-module, and $\text{Inn}_{(\sigma, \tau)}(\mathcal{A})$ becomes its submodule. If $1$ is the unity in $\mathcal{A}$ and $D$ is a $(\sigma, \tau)$-derivation of $\mathcal{A}$, then $D(1) = 0$. 
\end{remark}

In view of the the above Definitions \ref{definition 2.1} and \ref{definition 2.2} and the Remark \ref{remark 2.3}, the outer $(\sigma, \tau)$-derivations are precisely the elements of the factor module $\text{Out}_{(\sigma, \tau)}(\mathcal{A}) = \frac{\mathcal{D}_{(\sigma, \tau)}(\mathcal{A})}{\text{Inn}_{(\sigma, \tau)}(\mathcal{A})}$. Also, note that the set $\mathcal{D}_{(\sigma, \tau)}(\mathcal{A}) \setminus \text{Inn}_{(\sigma, \tau)}(\mathcal{A})$ is the set of all non-inner $(\sigma, \tau)$-derivations of $\mathcal{A}$.

The following lemma establishes a connection between our notions of outer $(\sigma, \tau)$-\\derivations and non-inner $(\sigma, \tau)$-derivations of $\mathcal{A}$.

\begin{lemma}\label{lemma 2.4}
Let $T = \{D_{i} \in \mathcal{D}_{(\sigma, \tau)}(\mathcal{A}) \mid i \in I\}$ ($I$ some indexing set) be a left transversal of $\text{Inn}_{(\sigma, \tau)}(\mathcal{A})$ in $\mathcal{D}_{(\sigma, \tau)}(\mathcal{A})$ with $0$ as the coset representative of the coset $\text{Inn}_{(\sigma, \tau)}(\mathcal{A})$. Then the non-inner $(\sigma, \tau)$-derivations of $\mathcal{A}$ correspond to the elements in the set $\bigcup_{D_{i} \in T \setminus \{0\}} (D_{i} + \text{Inn}_{(\sigma, \tau)}(\mathcal{A}))$. More precisely, $\mathcal{D}_{(\sigma, \tau)}(\mathcal{A}) \setminus \text{Inn}_{(\sigma, \tau)}(\mathcal{A}) = \bigcup_{D_{i} \in T \setminus \{0\}} (D_{i} + \text{Inn}_{(\sigma, \tau)}(\mathcal{A}))$.
\end{lemma}
\begin{proof}
Let $D \in \mathcal{D}_{(\sigma, \tau)}(\mathcal{A}) \setminus \text{Inn}_{(\sigma, \tau)}(\mathcal{A})$. Then $D \in \mathcal{D}_{(\sigma, \tau)}(\mathcal{A})$ but $D \notin \text{Inn}_{(\sigma, \tau)}(\mathcal{A})$. This implies that $D + \text{Inn}_{(\sigma, \tau)}(\mathcal{A}) \neq \text{Inn}_{(\sigma, \tau)}(\mathcal{A})$ so that $D + \text{Inn}_{(\sigma, \tau)}(\mathcal{A}) = D_{i} + \text{Inn}_{(\sigma, \tau)}(\mathcal{A})$ for some $D_{i} \in T \setminus \{0\}$. Therefore, $D - D_{i} \in \text{Inn}_{(\sigma, \tau)}(\mathcal{A})$ so that $D - D_{i} = D_{0}$ for some $D_{0} \in \text{Inn}_{(\sigma, \tau)}(\mathcal{A})$. Finally, we get that $D = D_{i} + D_{0}$ so that $D \in D_{i} + \text{Inn}_{(\sigma, \tau)}(\mathcal{A})$. Hence, $D \in \bigcup_{D_{i} \in T \setminus \{0\}} (D_{i} + \text{Inn}_{(\sigma, \tau)}(\mathcal{A}))$. Therefore, $\mathcal{D}_{(\sigma, \tau)}(\mathcal{A}) \setminus \text{Inn}_{(\sigma, \tau)}(\mathcal{A}) \subseteq \bigcup_{D_{i} \in T \setminus \{0\}} (D_{i} + \text{Inn}_{(\sigma, \tau)}(\mathcal{A}))$.

Conversely, let $D \in \bigcup_{D_{i} \in T \setminus \{0\}} (D_{i} + \text{Inn}_{(\sigma, \tau)}(\mathcal{A}))$. So $D \in D_{i} + \text{Inn}_{(\sigma, \tau)}(\mathcal{A})$ for some $D_{i} \in T \setminus \{0\}$. So $D + \text{Inn}_{(\sigma, \tau)}(\mathcal{A}) = D_{i} + \text{Inn}_{(\sigma, \tau)}(\mathcal{A})$. But since $D_{i} + \text{Inn}_{(\sigma, \tau)}(\mathcal{A}) \neq \text{Inn}_{(\sigma, \tau)}(\mathcal{A})$, so $D + \text{Inn}_{(\sigma, \tau)}(\mathcal{A}) \neq \text{Inn}_{(\sigma, \tau)}(\mathcal{A})$ so that $D \notin \text{Inn}_{(\sigma, \tau)}(\mathcal{A})$, that is, $D \in \mathcal{D}_{(\sigma, \tau)}(\mathcal{A}) \setminus \text{Inn}_{(\sigma, \tau)}(\mathcal{A})$. Therefore, $\bigcup_{D_{i} \in T \setminus \{0\}} (D_{i} + \text{Inn}_{(\sigma, \tau)}(\mathcal{A})) \subseteq \mathcal{D}_{(\sigma, \tau)}(\mathcal{A}) \setminus \text{Inn}_{(\sigma, \tau)}(\mathcal{A})$.
\end{proof}

\begin{remark}\label{remark 2.5}
In view of Lemma \ref{lemma 2.4}, $D$ is a non-inner derivation of $\mathcal{A}$ if and only if $D + \text{Inn}_{(\sigma, \tau)}(\mathcal{A})$ is a non-zero element of $\text{Out}_{(\sigma, \tau)}(\mathcal{A}) = \frac{\mathcal{D}_{(\sigma, \tau)}(\mathcal{A})}{\text{Inn}_{(\sigma, \tau)}(\mathcal{A})}$. In other words, $D$ is a non-inner derivation of $\mathcal{A}$ if and only if $D + \text{Inn}_{(\sigma, \tau)}(\mathcal{A})$ is a non-trivial outer derivation of $\mathcal{A}$. Therefore, studying the non-trivial outer derivations of $\mathcal{A}$ is equivalent to studying the non-inner derivations of $\mathcal{A}$ (see {\cite[Chapter 11]{pierce}} for details).
\end{remark} 

Many authors have defined an outer derivation of a ring $R$ (or algebra $\mathcal{A}$) to be a non-inner derivation of $R$ (or $\mathcal{A}$). Some references in this regard are \cite{batty1978derivations, chuang2005identities, dhara2022note, eroǧlu2017images, hall1972derivations, miles1964derivations, prajapati2022b, sakai1966derivations, weisfeld1960derivations}. In \cite{arutyunov2019smooth} and \cite{mishchenko2020description}, the authors initially describe the set $\mathcal{D}_{(\sigma, \tau)}(\mathcal{A}) \setminus \text{Inn}_{(\sigma, \tau)}(\mathcal{A})$ as the set of outer derivations but then argue that it is more natural to call the quotient module $\text{Out}_{(\sigma, \tau)}(\mathcal{A}) = \frac{\mathcal{D}_{(\sigma, \tau)}(\mathcal{A})}{\text{Inn}_{(\sigma, \tau)}(\mathcal{A})}$ as the set of outer derivations because this module can be interpreted as $1^{\text{st}}$ Hochschild cohomology module of $\mathcal{A}$ with coefficients in $\mathcal{A}$.


\begin{definition}\label{definition 2.6}
$\mathcal{A}$ is said to be $(\sigma, \tau)$-differentially trivial if $\mathcal{A}$ has only zero $(\sigma, \tau)$-derivation, that is, $\mathcal{D}_{(\sigma, \tau)}(\mathcal{A}) = \{0\}$.
\end{definition}

\begin{definition}\label{definition 2.7}
If $R$ is a ring and $G$ is a group, then the group ring of $G$ over $R$ is defined as the set $$RG = \{\sum_{g \in G} a_{g} g \mid a_{g} \in R, \forall g \in G \hspace{0.2cm} \text{and} \hspace{0.2cm} |\text{supp}(\alpha)| < \infty \},$$ where for $\alpha = \sum_{g \in G} a_{g} g$, $\text{supp}(\alpha)$ denotes the support of $\alpha$ that consists of elements from $G$ that appear in the expression of $\alpha$. $RG$ is a ring concerning the componentwise addition and multiplication defined, respectively, by: For $\alpha = \sum_{g \in G} a_{g} g$, $\beta = \sum_{g \in G} b_{g} g$ in $RG$, $$(\sum_{g \in G} a_{g} g ) + (\sum_{g \in G} b_{g} g) = \sum_{g \in G}(a_{g} + b_{g}) g \hspace{0.2cm} \text{and} \hspace{0.2cm} \alpha \beta = \sum_{g, h \in G} a_{g} b_{h} gh.$$ If $R$ is a commutative unital ring, then the group ring $RG$ is usually called a group algebra. If the ring $R$ is commutative, having unity $1$, and the group $G$ is abelian, having identity $e$, then $RG$ becomes a commutative unital algebra over $R$ with identity $1 = 1e$. We adopt the convention that empty sums are $0$, and empty products are $1$.
\end{definition}

\begin{definition}\label{definition 2.8}
The characteristic of a ring $R$ is the least positive integer $m$ such that $ma = 0$ for all $a \in R$. If no such integer exists, we say that $R$ has characteristic $0$. 
\end{definition}

\begin{definition}\label{definition 2.9}
Let $R$ be a unital ring and $G$ be a group. The elements of the form $rg$, where $r \in \mathcal{U}(R)$, the multiplicative group of units of $R$, and $g \in G$, which have a multiplicative inverse in $RG$, namely, $r^{-1}g^{-1}$, are called the trivial units of $RG$.
\end{definition}

\begin{notation}\label{notation 2.10}
The following notations are used throughout: For any $m, n \in \mathbb{N}$, $\text{gcd}(m,n)$ denotes the greatest positive common divisor of $m$ and $n$; we adopt the convention that $\text{gcd}(0,m) = m$ for any integer $m$.
$\text{char}(R)$ denotes the characteristic of ring $R$.  $\mathbb{N}$ denotes the set of natural numbers. $\mathbb{N}_{0} = \mathbb{N} \cup \{0\}$. $R$ denotes a commutative ring with unity. $\text{det}(A)$ denotes the determinant of a square matrix $A$. $A^{T}$ denotes the transpose of matrix $A$. For a group $G$, $|g|$ denotes the order of an element $g \in G$.
\end{notation}

Unless specifically mentioned, in the coming results, $R$ denotes a commutative ring with unity, $G = \{g_{1}, g_{2}, ..., g_{n}\}$ a finite group of order $n$, and ($\sigma, \tau$) a pair of unital (that is, $\sigma(1)=\tau(1)=1$) $R$-algebra endomorphisms of $RG$. First, we have the following lemma:
\begin{lemma}\th\label{lemma 2.11}
An $R$-linear map $D:RG \rightarrow RG$ is a $(\sigma, \tau)$-derivation if and only if $$D(g_{i}g_{j}) = D(g_{i}) \tau(g_{j}) + \sigma(g_{i}) D(g_{j})$$ for every $i, j \in \{1, 2, ..., n\}$.
\end{lemma}

For $k \in \mathbb{N}_{0}$, define the set $S_{k} = \{(i,j) \in \mathbb{N}_{0} \times \mathbb{N}_{0} \mid i+j = k\}$. For each $k \in \mathbb{N}_{0}$, $|S_{k}| = k+1$. Also, note that the sets $S_{k} \text{'s}$ are pairwise disjoint.

The following lemma gives a necessary condition for an $R$-linear map $D$ on $RG$ to be a $(\sigma, \tau)$-derivation of $RG$.
\begin{lemma}\th\label{lemma 2.12}
Let $D: RG \rightarrow RG$ be an $R$-linear map. If $D$ is a $(\sigma, \tau)$-derivation of $RG$, then $D(\alpha^{k}) = \left(\sum_{(i,j) \in S_{k-1}}  \sigma(\alpha^{i}) \tau(\alpha^{j})\right)D(\alpha)$ for every $\alpha \in RG$ and for every $k \in \mathbb{N}$.
\end{lemma}

\begin{proof} We use induction on $k$. For $k = 1$, the equality trivially holds. Let the result hold for $k=n$, that is, $D(\alpha^{n}) = \left(\sum_{(i,j) \in S_{n-1}}  \sigma(\alpha^{i}) \tau(\alpha^{j})\right)D(\alpha)$. Then

\begin{eqnarray*}
D(\alpha^{n+1}) = D(\alpha^{n} \alpha) & = & D(\alpha^{n}) \tau(\alpha) + \sigma(\alpha^{n}) D(\alpha) \\ & = & \left(\sum_{i+j = n-1} \sigma(\alpha^{i}) \tau(\alpha^{j+1}) + \sigma(\alpha^{n}) \tau(\alpha^{0})\right) D(\alpha) \\ & = & \left(\sum_{i+j = n} \sigma(\alpha^{i}) \tau(\alpha^{j})\right)D(\alpha) = \left(\sum_{(i,j) \in S_{n}}  \sigma(\alpha^{i}) \tau(\alpha^{j})\right)D(\alpha).
\end{eqnarray*}

Induction is complete, and so is proof.
\end{proof}

\begin{corollary}\th\label{corollary 2.13}
Let $G = \langle x \mid x^{n} = 1\rangle$ be a finite cyclic group of order $n$ and $D: RG \rightarrow RG$ be an $R$-linear map. If $D$ is a $(\sigma, \tau)$-derivation of $RG$, then \begin{equation}\label{eq 2.1}
D(x^{k}) = \left(\sum_{(i,j) \in S_{k-1}}  \sigma(x^{i}) \tau(x^{j})\right)D(x)\end{equation} for all $k \in \{1, 2, ..., n-1\}$.
\end{corollary}

\begin{remark}\th\label{remark 2.14} The converse of \th\ref{corollary 2.13} is not true in general. As a simple example, consider $G = C_{n} = \langle x \mid x^{n} = 1\rangle$, a cyclic group of order $n \geqslant 2$. Let $\sigma$ and $\tau$ be defined by $\sigma(x) = 1$ and $\tau(x) = x$. Let $D:RC_{n} \rightarrow RC_{n}$ be an $R$-linear map with $D(1) = 0$ and satisfying (\ref{eq 2.1}). Then $D(x^{n}) = D(1) = 0$ and $\left(\sum_{(i,j) \in S_{n-1}}  \sigma(x^{i}) \tau(x^{j})\right)D(x) = \left(1 + x + ... + x^{n-1}\right)D(x)$.
Note that $1 + x + ... + x^{n-1} \neq 0$ since $1, x, ..., x^{n-1}$ are linearly independent over $R$ being basis elements of $RC_{n}$. So if $D(x) \in RC_{n}$ is such that $\left(1 + x + ... + x^{n-1}\right)D(x) \neq 0$ (for example, $D(x) = x$), then $$D(x^{n}) \neq \left(\sum_{(i,j) \in S_{n-1}}  \sigma(x^{i}) \tau(x^{j})\right)D(x).$$ Hence, $D$ cannot be a $(\sigma, \tau)$-derivation of $RC_{n}$.
\end{remark}

In the following results, we assume that $\sigma$ and $\tau$ are distinct.

The following lemma gives a necessary and sufficient condition for a $(\sigma, \tau)$-derivation $D$ of $RG$ to be inner.
\begin{lemma}\th\label{lemma 2.15}
A $(\sigma, \tau)$-derivation $D:RG \rightarrow RG$ is inner if and only if there exists some $\beta \in RG$ such that $D(g_{i}) = \beta (\tau - \sigma)(g_{i})$ for all $i \in \{1, 2, ..., n\}$.
\end{lemma}
\begin{proof}
For $\alpha = \sum_{i=1}^{n} a_{i} g_{i} \in RG$, \begin{eqnarray*}D(\alpha) = \sum_{i=1}^{n} a_{i} D(g_{i}) & = & \sum_{i=1}^{n} a_{i} \beta (\tau - \sigma)(g_{i}) = \beta \left(\sum_{i=1}^{n} a_{i} \tau(g_{i}) - \sum_{i=1}^{n} a_{i} \sigma(g_{i})\right) \\ & = & \beta \left(\tau\left(\sum_{i=1}^{n} a_{i} g_{i}\right) - \sigma\left(\sum_{i=1}^{n} a_{i} g_{i}\right)\right) = \beta (\tau - \sigma)(\alpha).\end{eqnarray*}
Since $\alpha \in \mathcal{A}$ is arbitrary, therefore, $D$ is inner. The forward part follows directly from the definition of an inner $(\sigma, \tau)$-derivation.
\end{proof}

\begin{lemma}\th\label{lemma 2.16}
Let $D: RG \rightarrow RG$ be a $(\sigma, \tau)$-derivation and $\alpha \in RG$. If there exists some $\beta \in RG$ such that $D(\alpha) = \beta (\tau - \sigma)(\alpha)$, then $D(\alpha^{n}) = \beta (\tau - \sigma)(\alpha^{n})$ for all $n \in \mathbb{N}$.
\end{lemma}

\begin{proof} Follows by induction on $n$.
\end{proof}

\begin{corollary}\th\label{corollary 2.17}
Let $G = \langle x \mid x^{n} = 1\rangle$ be a finite cyclic group of order $n$ and $D:RG \rightarrow RG$ be a $(\sigma, \tau)$-derivation. Then $D$ is inner if and only if there exists some $\beta \in RG$ such that $D(x) = \beta (\tau - \sigma)(x)$.
\end{corollary}
\begin{proof}
The forward part follows trivially by definition. For the converse, note that $D(1) = D(1) \tau(1) + \sigma(1) D(1) = D(1) + D(1)$, so that $D(1) = 0$. Also, $\beta(\tau - \sigma)(1) = 0$ so that $D(1) = \beta(\tau - \sigma)(1)$. Further, in view of \th\ref{lemma 2.16}, $D(x^{k}) = \beta (\tau - \sigma)(x^{k})$ for all $k \in \{1, ..., n-1\}$. Now, using \th\ref{lemma 2.15}, $D$ is inner.
\end{proof}

Our proof techniques will be purely algebraic, and we will significantly use linear algebra over commutative rings. We refer the reader to \cite{McDonald1984} and \cite{Payne2009} for the theory of linear algebra over commutative rings.

\section{Inner $(\sigma, \tau)$-Derivations of Cyclic Group Rings}\label{section 4}
Throughout the section, unless otherwise stated, we assume that $R$ is a commutative ring with unity, $G = \langle x \mid x^{n} = 1 \rangle$ is a cyclic group of order $n$, and $(\sigma, \tau)$ is a pair of distinct $R$-algebra endomorphisms of $RG$, which are $R$-linear extensions of group endomorphisms of $G$. This section determines a necessary and sufficient condition for a $(\sigma, \tau)$-derivation $D$ of $RG$ to be inner. 

Since $\sigma, \tau:G \rightarrow G$ are distinct group homomorphisms, therefore, $\sigma(x) = x^{u}$ and $\tau(x) = x^{u+v}$ for some $u, v \in \{0, 1, ..., n-1\}$ with $v \neq 0$.

In this section, we give a characterization for a $(\sigma, \tau)$-derivation $D$ of $RG$ to be inner in the main theorem, namely, \th\ref{theorem 4.1} of this section. We prove this main \th\ref{theorem 4.1} through a sequence of lemmas.

\begin{lemma}\th\label{lemma 4.2}
Suppose that $\beta = \sum_{i=0}^{n-1}b_{i}x^{i}$. Then solving the equation $\beta (\sigma - \tau)(x) = D(x)$ for $\beta$ in $RG$ is the same as solving the matrix equation $AB = C$ for $B$, where $A$ is an $n \times n$ matrix over $R$ whose $(i,j)^{\text{th}}$-entry ($i, j \in \{0, 1, ..., n-1\}$) is given by \begin{equation}\label{eq 4.2}
a_{i,j} = \begin{cases}
1 & \text{if $j \equiv (i-u) (\text{mod} ~ n)$} \\
-1 & \text{if $j \equiv (i-u-v) (\text{mod} ~ n)$} \\
0 & \text{otherwise}
\end{cases},\end{equation} $B = \begin{pmatrix}
b_{0} & b_{1} & \cdots & b_{n-1}
\end{pmatrix}^{T}$ and $C = \begin{pmatrix}
c_{0} & c_{1} & \cdots & c_{n-1}\end{pmatrix}^{T}$.
\end{lemma}
\begin{proof}
As $\beta \in RG$, so $\beta = \sum_{i=0}^{n-1}b_{i}x^{i}$ for some $b_{i} \in R$ ($0 \leq i \leq n-1$). We need to solve the matrix equation $AB = C$ for the unknown $B$ consisting of the $n$ unknown distinct variables $b_{0}, b_{1}, ..., b_{n-1}$. Also assume that $b_{i} = b_{j}$ whenever $i \equiv j$ (mod n). Now, using our assumption on $b_{i}$'s, we get that 

\begin{eqnarray*}
\beta (\sigma - \tau)(x) & = & \left(\sum_{i=0}^{n-1}b_{i}x^{i}\right)(x^{u} - x^{u+v}) = \sum_{j=u}^{u+n-1}b_{j-u}x^{j} - \sum_{k=u+v}^{u+v+n-1} b_{k-u-v}x^{k} \\ & = & \sum_{i=0}^{n-1} b_{i-u} x^{i} - \sum_{i=0}^{n-1} b_{i-u-v}x^{i} = \sum_{i=0}^{n-1} (b_{i-u} - b_{i-u-v})x^{i}.
\end{eqnarray*}

Finally, using the fact that $G = \{x^{i} \mid 0 \leq i \leq n-1\}$ is an $R$-basis of $RG$, we get that $\beta(\sigma - \tau)(x) = D(x)$ if and only if $b_{i-u} - b_{i-u-v} = c_{i}$ for all $i \in \{0, 1, ..., n-1\}$.
\end{proof}

\begin{lemma}\th\label{lemma 4.3}
The augmented matrix $\left(\begin{array}{c|c}
A & C
\end{array}\right)$ is row equivalent to the matrix $\left(\begin{array}{c|c}
A' & C' \\
0 & C''
\end{array}\right)$ where $A'$ is the $(n-d) \times n$ matrix whose $n-d$ rows are the same as the first $n-d$ rows of $A$, $C'$ is the $(n-d) \times 1$ matrix given by $C' = \begin{pmatrix}
c_{0} &
c_{1} &
\cdots &
c_{n-d-1}
\end{pmatrix}^{T}$ and $C''$ is the $d \times 1$ matrix given by $C'' = \begin{pmatrix}
\sum_{j=0}^{m-1} c_{dj} &
\sum_{j=0}^{m-1} c_{1+dj} &
\cdots &
\sum_{j=0}^{m-1} c_{(d-1)+dj}
\end{pmatrix}^{T}$.
\end{lemma}
\begin{proof}
Let $R_{0}, R_{1},..., R_{n-1}$ denote the $n$ rows of the $n \times n$ matrix $A$ in the given order. Then, for each $i \in \{0, 1, ..., n-1\}$, $R_{i}B = b_{i-u} - b_{i-u-v}$. It is enough to prove that for each $i \in \{0, 1, ..., d-1\}$, $\sum_{j=0}^{m-1}R_{i+jd}B = 0$. Now, for $i \in \{0, 1, ..., d-1\}$, \begin{eqnarray*}
\sum_{j=0}^{m-1}R_{i+jd}B & = & \sum_{j=0}^{m-1} (b_{i+jd-u} - b_{i+jd-u-v}) = \sum_{j=0}^{m-1} b_{i+jd-u} - \sum_{j=0}^{m-1} b_{i+jd-u-v}.
\end{eqnarray*}

Let $j, j' \in \{0, 1, ..., m-1\}$. Also, put $m_{0} = \frac{v}{d}$. Note that $1 \leq m_{0} \leq m$. Then \begin{equation*}
\begin{aligned}
& b_{i+jd-u} = b_{i+j'd-u-v} 
\\ & \Leftrightarrow (i+jd-u) \equiv (i+j'd-u-v)(\text{mod} ~ n) 
\\ & \Leftrightarrow jd \equiv (j'd-v) (\text{mod} ~ n)
\\ & \Leftrightarrow j'd \equiv (jd + m_{0}d) (\text{mod} ~ n) 
\\ & \Leftrightarrow j' \equiv (j + m_{0}) (\text{mod} ~ m).
\end{aligned}
\end{equation*} 

Therefore, for each $j \in \{0, 1, ..., m-1\}$, there exists a unique $j' \in \{0, 1, ... m-1\}$ satisfying $j' \equiv (j + m_{0}) (\text{mod} ~ m)$ such that $b_{i+jd-u} = b_{i+j'd-u-v}$. Therefore, $\sum_{j=0}^{m-1}R_{i+jd}B = 0$.
\end{proof}

It is important to note that to obtain the above row equivalent matrix form of $A$, we used only one type of elementary row operation: the addition of rows.

\begin{lemma}\th\label{lemma 4.4}
The first $n-d$ rows $R_{0}, R_{1}, ..., R_{n-d-1}$ of $A$ are linearly independent over $R$. 
\end{lemma}
\begin{proof}
$R^{n} = \underbrace{R \times ... \times R}_{n-\text{times}}$ is a ring concerning the usual operations of component-wise addition and multiplication, and hence naturally is an $R$-module. The rows $R_{0}, R_{1}, ..., R_{n-d-1}$ of $A$ can be considered elements of $R^{n}$. To show that the rows $R_{0}, R_{1}, ..., R_{n-d-1}$ are linearly independent over $R$, it is to be proved that if $\alpha_{i} \in R$ ($0 \leq i \leq n-d-1$) such that $\sum_{i=0}^{n-d-1}\alpha_{i}R_{i} = 0$, then $\alpha_{i} = 0$ for all $i \in \{0, 1, ..., n-d-1\}$. So, let us assume that $\sum_{i=0}^{n-d-1}\alpha_{i}R_{i} = 0$ for some $\alpha_{i} \in R$ ($0 \leq i \leq n-d-1$).

Note that $\sum_{i=0}^{n-d-1}\alpha_{i}R_{i} = 0$ is equivalent to $n$ equations in $n-d$ variables $\alpha_{i}$ ($0 \leq i \leq n-d-1$): \begin{equation}\label{eq 4.3}
\sum_{i=0}^{n-d-1} \alpha_{i} a_{i,j} = 0, ~~~~~ j \in \{0, 1, ..., n-1\}.
\end{equation}

Let $i \in \{0, 1, ..., n-d-1\}$ be arbitrary but fixed. We show that $\alpha_{i} = 0$. 

It is an easy observation from (\ref{eq 4.2}) that each row and each column of $A$ contains precisely one 1 and one -1. So, there exist unique $j_{1}, j_{2} \in \{0,1, ..., n-1\}$, $j_{1} \neq j_{2}$ such that $a_{i,j_{1}} = 1$, $a_{i,j_{2}} = -1$, and $a_{i,j} = 0$ for all $j \in \{0, 1, ..., n-1\}$ with $j \neq j_{1}, j_{2}$.

We have some possible cases.

\textbf{Case 1:} There exists some $i_{1} \in \{n-d, ..., n-1\}$ such that $a_{i_{1},j_{1}} = -1$.

Then, using the fact that each column (in particular, the $j_{1}^{\text{th}}$ column) of $A$ contains exactly two non-zero entries, namely, $1$ and $-1$, we get from (\ref{eq 4.3}) that $\alpha_{i} = 0$.

\textbf{Case 2:} There exists some $i_{2} \in \{n-d, ..., n-1\}$ such that $a_{i_{2},j_{2}} = 1$.

Again, using the fact that the $j_{2}^{\text{th}}$ column of $A$ contains only two non-zero entries, namely, $1$ and $-1$, we will get from (\ref{eq 4.3}) that $-\alpha_{i} = 0$ so that $\alpha_{i} = 0$.

\textbf{Case 3:} There do not exist any $i_{1}, i_{2} \in \{n-d, ..., n-1\}$ such that $a_{i_{1},j_{1}} = -1$ and $a_{i_{2},j_{2}} = 1$.

Then, because each column (in particular, the $j_{1}^{\text{th}}$ and $j_{2}^{\text{th}}$ columns) of $A$ contains only two non-zero entries, $1$ and $-1$, we will get that there will exist some $i', i'' \in \{0, 1, ..., n-d-1\}$, $i', i'' \neq i$ such that $a_{i', j_{1}} = -1$ and $a_{i'',j_{2}} = 1$. Therefore, by (\ref{eq 4.3}), we have two equations: $\alpha_{i} - \alpha_{i'} = 0$ and $- \alpha_{i} + \alpha_{i''} = 0$, so that $\alpha_{i} = \alpha_{i'} = \alpha_{i''}$.

\textbf{Claim:} $i'$ and $i''$ are distinct.

If possible, suppose that $i' = i''$. Then $R_{i} + R_{i'} = 0$. Since $a_{i,j_{1}} = 1$ and $a_{i', j_{1}} = -1$, therefore, by (\ref{eq 4.2}), $j_{1} \equiv (i-u) (\text{mod} ~ n)$ and $j_{1} \equiv (i' - u - v) (\text{mod} ~ n)$ so that $(i-u) \equiv (i'-u-v) (\text{mod} ~ n)$, which holds if and only if $v \equiv (i'-i) (\text{mod} ~ n)$. Further, since $a_{i,j_{2}} = -1$ and $a_{i', j_{2}} = 1$, therefore, by (\ref{eq 4.2}), $j_{2} \equiv (i-u-v) (\text{mod} ~ n)$ and $j_{2} \equiv (i' - u) (\text{mod} ~ n)$ so that $(i-u-v) \equiv (i'-u) (\text{mod} ~ n)$, which holds if and only if $v \equiv (i-i') (\text{mod} ~ n)$. So, we get that $v \equiv -v (\text{mod} ~ n)$, which implies that $2v \equiv 0 (\text{mod} ~ n)$, or equivalently, $v \equiv 0 (\text{mod} ~ \text{gcd}(n,2))$. If $n$ is odd, this gives us that $v \equiv 0 (\text{mod} ~ n)$ so that $v=0$, which is not possible because of our assumption that $\sigma$ and $\tau$ are distinct so that $v \neq 0$. So, we arrive at a contradiction. If $n$ is even, we get that $v \equiv 0 (\text{mod} ~ \frac{n}{2})$ so that $v = \frac{n}{2}$ (because $v \in \{1, ..., n-1\}$). This means that $d = \text{gcd}(v,n) = \text{gcd}(\frac{n}{2},n) = \frac{n}{2} = v$, so that $n-d = \frac{n}{2}$. Therefore, from the above congruences, we will get that $i' \equiv (i + \frac{n}{2})(\text{mod} ~ n)$. But this congruence has no solution for $i, i' \in \{0, 1, ..., n-d-1 = \frac{n}{2} - 1\}$. So, we again arrive at a contradiction. Therefore, our assumption regarding the equality of $i'$ and $i''$ must be wrong. In other words, $i'$ and $i''$ are distinct. Hence, the claim is proven.

Now, consider $i'$. Let $j' \in \{0, 1, ..., n-1\}$ be such that $a_{i',j'} = 1$. Note that $j' \neq j_{1}, j_{2}$. Check if there exists some $i^{(3)} \in \{n-d, ..., n-1\}$ such that $a_{i^{(3)}, j'} = -1$. If yes, then using the idea in Case 1, we will get that $\alpha_{i'} = 0$, so that $\alpha_{i} = 0$. If no such $i^{(3)} \in \{n-d, ..., n-1\}$ exists, then consider $i''$. If $j'' \in \{0, 1, ..., n-1\}$ is such that $a_{i'',j''} = -1$ (clearly $j'' \neq j_{1}, j_{2}$), then check if there exists some $i^{(4)} \in \{n-d, ..., n-1\}$ such that $a_{i^{(4)},j''} = 1$. If yes, then using the idea in Case 2, we will get that $\alpha_{i''} = 0$, so that $\alpha_{i} = 0$. If not, consider either $i'$ or $i''$ and do the following: Suppose we consider $i'$. Since there does not exist any $i^{(3)} \in \{n-d, ..., n-1\}$ such that $a_{i^{(3)}, j'} = -1$, there will exist some $i^{(3)} \in \{0, 1, ..., n-d-1\}$ such that $a_{i^{(3)}, j'} = -1$. In view of (\ref{eq 4.3}), this gives us that $\alpha_{i'} = \alpha_{i^{(3)}}$. Now, consider the $(i^{(3)})^{\text{th}}$ row. Let $j^{(3)} \in \{0, 1, ..., n-1\}$ be such that $a_{i^{(3)}, j^{(3)}} = 1$. Note that $j^{(3)} \neq j' \neq j_{1} \neq j_{2}$. Now, check if there exists some $i^{(4)} \in \{n-d, ..., n-1\}$ such that $a_{i^{(4)}, j^{(3)}} = -1$. If yes, then from (\ref{eq 4.3}), we will get that $\alpha_{i^{(3)}} = 0$ so that $\alpha_{i} = 0$. If not, then there will exist some $i^{(4)} \in \{0, 1, ..., n-d-1\}$ such that $a_{i^{(4)},j^{(3)}} = -1$. Now, we consider the $(i^{(4)})^{\text{th}}$ row and repeat the same process as we did for $i^{(3)}$. Now, keep repeating the same process until we get a row $i^{(r)} \in \{0, 1, ..., n-d-1\}$ and a column $j^{(r)} \in \{0, 1, ..., n-1\}$ such that there will exist some $i^{(r+1)} \in \{n-d, ..., n-1\}$ with $a^{i^{(r+1)},j^{(r)}} = -1$ thus giving us $\alpha_{i^{(r)}} = 0$ so that $\alpha_{i} = 0$.
\end{proof}

\begin{lemma}\th\label{lemma 4.5}
The coefficient matrix $A$ defined in (\ref{eq 4.2}) is a non-invertible (that is, $\text{det}(A) = 0$) circulant matrix of row rank $n-d$. 
\end{lemma}
\begin{proof}
The first row $R_{0}$ of $A$ corresponds to $i=0$. Define $\theta:R^{n} \rightarrow R^{n}$ by $\theta(a_{0}, a_{1}, ..., a_{n-1}) = (a_{n-1}, a_{0}, a_{1}, ..., a_{n-2})$ for all $(a_{0}, a_{1}, ..., a_{n-1}) \in R^{n}$. Then $\theta$ is an $R$-module homomorphism, and $\theta$ shifts each entry of an element in $R^{n}$ to the right by one step. Observe that $\theta(R_{0}) = R_{1}$ and $\theta(R_{1}) = R_{2}$, so that $\theta^{2}(R_{0}) = R_{2}$. Likewise, $\theta^{i}(R_{0}) = R_{i}$ for all $i \in \{0, 1, ..., n-1\}$. In other words, $\theta^{i}$ shifts the vector $R_{0}$ to the right by $i$ number of steps. Therefore, A is an $n \times n$ circulant matrix. 

Since $A$ is row equivalent to a matrix whose last $d$ number of rows is zero, and from the observation made just after the proof of \th\ref{lemma 4.3}, we know that this row equivalent matrix form was obtained using just one type of elementary row operation, namely, the addition of rows, therefore, the determinant of $A$ is equal to the determinant of the obtained row equivalent matrix form {\cite[Lemma 6.1.7 (b)]{Payne2009}}. However, since the row equivalent matrix form contains at least one zero row, its determinant is zero {\cite[Lemma 6.1.7 (a)]{Payne2009}}. Hence, the determinant of $A$ is zero. Hence, $A$ is a non-invertible matrix (using {\cite[Theorem 6.3.7]{Payne2009}}).

By \th\ref{lemma 4.3}, the row rank of $A$ is at most $n-d$. Further, by \th\ref{lemma 4.4}, the row rank of $A$ is precisely $n-d$ as its first $n-d$ rows are linearly independent over $R$. 
\end{proof}

We now prove the following \th\ref{theorem 4.1} using the above sequence of lemmas.
\begin{theorem}\th\label{theorem 4.1}
Let $R$ be a commutative unital ring, $G = \langle x \mid x^{n} = 1\rangle$ be a cyclic group of order $n$, and $(\sigma, \tau)$ be a pair of $R$-algebra endomorphisms of $RG$ which are $R$-linear extensions of the group endomorphisms of $G$ (so that $\sigma(x) = x^{u}$ and $\tau(x) = x^{u+v}$ for some $u, v \in \{0, 1, ..., n-1\}$). Let $D:RG \rightarrow RG$ be a $(\sigma, \tau)$-derivation. Suppose that $D(x) = \sum_{i=0}^{n-1}c_{i}x^{i}$, $d = \text{gcd}(v,n)$, and $m = \frac{n}{d}$. Then $D$ is inner if and only if the following $d$ equations hold simultaneously: \begin{equation}\label{eq 4.1}
\sum_{j=0}^{m - 1} c_{i+jd} = 0, ~~~~~ i \in \{0, 1, ..., d-1\}.
\end{equation}
\end{theorem}
\begin{proof}
By \th\ref{corollary 2.17}, $D$ is inner if and only if there exists some $\beta \in RG$ such that $\beta(\sigma - \tau)(x) = D(x)$. So, to determine if $D$ is inner, we need to find if the equation $\beta(\sigma - \tau)(x)$ is solvable for $\beta$ in $RG$. 

By \th\ref{lemma 4.2}, solving the equation $\beta(\sigma - \tau)(x)$ for $\beta$ in $RG$ is the same as solving the matrix equation $AB=C$ for unknown $B$ with entries in $R$, where $A$, $B$, and $C$ are as in \th\ref{lemma 4.2}.

Further by \th\ref{lemma 4.3}, the augmented matrix $(A \mid C)$ is row equivalent to a matrix of the form $$\left(\begin{array}{c|c}
A' & C' \\
0 & C''
\end{array}\right) = \left(\begin{array}{ccc|c}
\cdots & R_{0} & \cdots & c_{0} \\
\cdots & R_{1} & \cdots & c_{1} \\
\vdots & \vdots & \vdots & \vdots \\
\cdots & R_{n-d-1} & \cdots & c_{n-d-1} \\
0 & \cdots & 0 & \sum_{j=0}^{m-1}c_{dj} \\
0 & \cdots & 0 & \sum_{j=0}^{m-1}c_{1+dj} \\
\vdots & \vdots & \vdots & \vdots \\
0 & \cdots & 0 & \sum_{j=0}^{m-1}c_{(d-1)+dj} \\
\end{array}\right),$$ where the matrices $A'$, $C'$, and $C''$ are as in \th\ref{lemma 4.3} and $R_{0}, R_{1}, ..., R_{n-d-1}$ are the first $n-d$ rows of the matrix $A$ in the given order. 

Further, since by \th\ref{lemma 4.4}, the first $n-d$ rows $R_{0}, R_{1}, ..., R_{n-d-1}$ of $A$ are linearly independent over $R$ and by \th\ref{lemma 4.5}, $A$ is a circulant matrix with each row and each column containing only two non-zero entries, namely, $1$ and $-1$, we can further row reduce the matrix $\left(\begin{array}{c|c}
A' & C' \\
0 & C''
\end{array}\right)$ into an $n \times n$ matrix of the form $\left(\begin{array}{cc|c}
I_{n-d} & A'' & C''' \\
0 & 0 & C''
\end{array}\right)$, where $I_{n-d}$ denotes the $(n-d) \times (n-d)$ identity matrix over $R$, $A''$ is an $(n-d) \times d$ matrix whose entries belong to $\{-1, 0, 1\}$, and $C'''$ is an $(n-d) \times 1$ matrix whose entries are linear combinations of $c_{i}$'s ($0 \leq i \leq n-d-1$). 

Combining all the above information, we get that $AB = C$ has a solution $B$ over $R$ if and only if the $d$ equations in  (\ref{eq 4.1}) hold simultaneously. 

Hence, the theorem is proved.
\end{proof}

By \th\ref{lemma 4.5}, the determinant of $A$ is $0$. Therefore, by corollary to {\cite[Theorem 1]{Ching1977}}, we have that if a solution of $AB = C$ exists, then it is never unique.

\section{Outer $(\sigma, \tau)$-Derivations of Cyclic Group Rings}\label{section 3}
Throughout this section, unless otherwise stated, $R$ denotes a commutative ring with unity, $G$ a cyclic group $G = \langle x \mid x^{n} = 1\rangle$ of order $n$, and $(\sigma, \tau)$ a pair of $R$-algebra endomorphisms of $RG$, which are $R$-linear extensions of group endomorphisms of $G$. In this section, we classify all $(\sigma, \tau)$-derivations of the cyclic group ring $RG$ (see \th\ref{lemma 3.1}). Consequently, we answer the twisted derivation problem for $RG$ in \th\ref{theorem 3.4}.

Since $\sigma, \tau: G \rightarrow G$ are group homomorphisms, we can write $\sigma(x) = x^{u}$ and $\tau(x) = x^{u+v}$ for some $u, v \in \{0, 1, ..., n-1\}$.

\begin{lemma}\th\label{lemma 3.1}
Let $R$ be a commutative unital ring, $G = \langle x \mid x^{n} = 1\rangle$ be a cyclic group of order $n$, and $(\sigma, \tau)$ be a pair of $R$-algebra endomorphisms of $RG$ which are $R$-linear extensions of the group endomorphisms of $G$ (so that $\sigma(x) = x^{u}$ and $\tau(x) = x^{u+v}$ for some $u, v \in \{0, 1, ..., n-1\}$). Let $D:RG \rightarrow RG$ be an $R$-linear map with $D(1) = 0$ and satisfying the $n-1$ equations in (\ref{eq 2.1}). Then, the following statements are equivalent:
\begin{enumerate}
\item[(i)] $D$ is a $(\sigma, \tau)$-derivation.
\item[(ii)] The following equation holds: \begin{equation}\label{eq 3.1}
\left(\sum_{i=0}^{n-1}x^{iv}\right)D(x) = 0.
\end{equation}
\item[(iii)] If $D(x) = \sum_{i=0}^{n-1}c_{i}x^{i}$ for some $c_{i} \in R$ ($0 \leq i \leq n-1$) and $d = \text{gcd}(v,n)$, then $d \left(\sum_{j=0}^{m-1} c_{i+jd} \right) = 0$ for all $i \in \{0, 1, ..., d-1\}$.
\end{enumerate}
\end{lemma}
\begin{proof}
First, we prove $(i) \Leftrightarrow (ii)$. In view of \th\ref{lemma 2.11}, $D$ will be a $(\sigma, \tau)$-derivation if and only if \begin{equation}\label{eq 3.2}D(x^{s+t}) = D(x^{s}) \tau(x^{t}) + \sigma(x^{s}) D(x^{t})\end{equation} for all $s, t \in \{0, 1, ..., n-1\}$.

The relations (\ref{eq 3.2}) hold trivially when at least one of $s$ or $t$ is $0$, using the fact that $\sigma(1) = \tau(1) = 1$ and $D(1) = 0$. So, now let $s, t \in \{1, 2, ..., n-1\}$. Using (\ref{eq 2.1}), we get:
\begin{eqnarray*}
D(x^{s}) \tau(x^{t}) = \left(\sum_{(i,j) \in S_{s-1}}  \sigma(x^{i}) \tau(x^{j})\right)D(x) \tau(x^{t}) & = & \left(\sum_{(i,j) \in S_{s-1}}  \sigma(x^{i}) \tau(x^{j+t})\right)D(x) \\ & = & \left(\sum_{i=0}^{s-1} \sigma(x^{i}) \tau(x^{s-1-i+t})\right)D(x);\end{eqnarray*}

\begin{eqnarray*}
\sigma(x^{s}) D(x^{t}) & = & \sigma(x^{s})\left(\sum_{(i,j) \in S_{t-1}} \sigma(x^{i}) \tau(x^{j})\right)D(x) = \left(\sum_{(i,j) \in S_{t-1}} \sigma(\zeta^{s+i}) \tau(x^{j})\right)D(x) \\ & = & \left(\sum_{i=0}^{t-1} \sigma(x^{s+i}) \tau(x^{t-1-i})\right)D(x) = \left(\sum_{i=s}^{s+t-1} \sigma(x^{i}) \tau(x^{t-1+s-i})\right)D(x).
\end{eqnarray*}
Therefore, \begin{equation}\label{eq 3.3}
D(x^{s}) \tau(x^{t}) + \sigma(x^{s}) D(x^{t}) = \left(\sum_{(i,j) \in S_{s+t-1}} \sigma(x^{i}) \tau(x^{j})\right)D(x).
\end{equation}

Since $s, t \in \{1, 2,..., n-1\}$, so $s+t \in \{2, 3, ..., 2(n-1)\}$. We split the proof into the following three cases:

\textbf{Case 1:} $s+t \leq n-1$.

Since $s+t \leqslant n-1$, by (\ref{eq 2.1}), $D(x^{s+t}) = \left(\sum_{(i,j) \in S_{s+t-1}} \sigma(x^{i}) \tau(x^{j})\right)D(x)$. Then, by (\ref{eq 3.3}), $D(x^{s}) \tau(x^{t}) + \sigma(x^{s}) D(x^{t}) = D(x^{s+t})$. Therefore, the relation (\ref{eq 3.2}) holds in this case.

\textbf{Case 2:} $s+t = n$.

Then $D(x^{s+t}) = D(x^{n}) = D(1) = 0$. And using (\ref{eq 3.3}),
\begin{eqnarray*}
D(x^{s}) \tau(x^{t}) + \sigma(x^{s}) D(x^{t}) = \left(\sum_{(i,j) \in S_{n-1}} x^{ui} x^{(u+v)j} \right)D(x) = x^{(n-1)u} \left(\sum_{j=0}^{n-1} x^{vj} \right)D(x).
\end{eqnarray*}

Therefore, $D(x^{s+t}) = D(x^{s}) \tau(x^{t}) + \sigma(x^{s}) D(x^{t})$ if and only if $\left(\sum_{j=0}^{n-1} x^{vj} \right)D(x) = 0$.

\textbf{Case 3:} $n < s+t \leq 2(n-1)$.

Then we can write $s+t = m + n$ for some $m \in \{1, ..., n-1\}$. Now by (\ref{eq 2.1}), 

\begin{eqnarray*}
D(x^{s+t}) & = & D(x^{m+n}) = D(x^{m}) = \left(\sum_{(i,j) \in S_{m-1}} \sigma(x^{i}) \tau(x^{j})\right)D(x) = x^{(m-1)u} \left(\sum_{j=0}^{m-1} x^{vj}\right)D(x).
\end{eqnarray*}

Further using (\ref{eq 3.3}), \begin{eqnarray*}
D(x^{s}) \tau(x^{t}) + \sigma(x^{s}) D(x^{t}) & = & x^{(m+n-1)u} \left(\sum_{j=0}^{m+n-1} x^{vj}\right)D(x) \\ & = & x^{(m-1)u} \left(\sum_{j=0}^{m-1} x^{vj}\right)D(x) + x^{(m-1)u + mv} \left(\sum_{j=0}^{n-1} x^{vj}\right)D(x).
\end{eqnarray*}


Therefore, $D(x^{s+t}) = D(x^{s}) \tau(x^{t}) + \sigma(x^{s}) D(x^{t})$ if and only if $\left(\sum_{j=0}^{n-1} x^{vj}\right)D(x) = 0$.

Combining all three cases, we can conclude that $D$ is a $(\sigma, \tau)$-derivation if and only if $\left(\sum_{j=0}^{n-1} x^{vj}\right)D(x) = 0$.

Now we prove $(ii) \Leftrightarrow (iii)$. Since $D(x) \in RG$, we can write $D(x) = \sum_{i=0}^{n-1} c_{i}x^{i}$ for some $c_{i} \in R$ ($0 \leq i \leq n-1$).

$v \in \{0, 1, ..., n-1\}$, so suppose $|x^{v}| = m$. Let $d = \frac{n}{m}$. Note that $d$ also equals $\text{gcd}(v,n)$. Then $\langle x^{v} \rangle = \langle x^{d} \rangle$. Put $c_{0} = c_{n}$. Now,

\begin{equation*}
\begin{aligned}
\left(\sum_{j=0}^{n-1} x^{vj} \right) D(x) & = \left(\sum_{j=1}^{n} x^{vj} \right) D(x) = d \left(\sum_{j=1}^{m} x^{vj} \right) D(x) = d \left( \sum_{i=1}^{md} c_{i} \left(\sum_{j=1}^{m} x^{dj + i} \right) \right)
\\ & = d [c_{1}(\sum_{j=1}^{m} x^{dj+1}) + c_{2} (\sum_{j=1}^{m} x^{dj+2}) + ... + c_{d} (\sum_{j=1}^{m} x^{dj+d})] + d[c_{d+1} (\sum_{j=1}^{m} x^{dj+d+1}) \\ &\quad + c_{d+2} (\sum_{j=1}^{m} x^{dj+d+2})  + ... + c_{d+d} (\sum_{j=1}^{m} x^{dj+d+d})] + ... \\ &\quad + d[c_{(m-1)d+1} (\sum_{j=1}^{m} x^{dj+(m-1)d+1}) + c_{(m-1)d+2}(\sum_{j=1}^{m} x^{dj+(m-1)d+2}) + ... \\ &\quad + c_{(m-1)d+d}(\sum_{j=1}^{m} x^{dj+(m-1)d+d})]
\\ & = d[c_{1}(\sum_{j=1}^{m} x^{dj+1}) + c_{2} (\sum_{j=1}^{m} x^{dj+2}) + ... + c_{d} (\sum_{j=1}^{m} x^{dj+d})] + d[c_{d+1} (\sum_{j=1}^{m} x^{dj+1}) \\ &\quad + c_{d+2} (\sum_{j=1}^{m} x^{dj+2}) + ... + c_{d+d} (\sum_{j=1}^{m} x^{dj+d})] + ... + d[c_{(m-1)d+1} (\sum_{j=1}^{m} x^{dj+1}) \\ &\quad + c_{(m-1)d+2}(\sum_{j=1}^{m} x^{dj+2}) + ... + c_{(m-1)d+d}(\sum_{j=1}^{m} x^{dj+d})]
\\ & = d(\sum_{i=0}^{m-1} c_{di+1})(\sum_{j=1}^{m}x^{dj+1}) + d(\sum_{i=0}^{m-1} c_{di+2})(\sum_{j=1}^{m}x^{dj+2}) + ... + \\ &\quad d(\sum_{i=0}^{m-1} c_{di+d})(\sum_{j=1}^{m}x^{dj+d}).
\end{aligned}
\end{equation*}

Therefore, $\left(\sum_{j=0}^{n-1} x^{vj} \right) D(x) = 0$ if and only if $d\left(\sum_{i=0}^{m-1} c_{di+j}\right) = 0$ for all $j \in \{1, 2, ..., d\}$ or, equivalently, $d\left(\sum_{i=0}^{m-1} c_{di+j}\right) = 0$ for all $j \in \{0, 1, ..., d-1\}$.

Therefore, $D$ is a $(\sigma, \tau)$-derivation if and only if $\sum_{i=0}^{m-1} c_{i+jd} = 0$ for all $j \in \{0, 1, ..., d-1\}$.

Hence, the theorem is proved.
\end{proof}

Some immediate corollaries to the above theorem can be stated as follows:
\begin{corollary}\th\label{corollary 3.2}
Let $D:RG \rightarrow RG$ be an $R$-linear map with $D(1) = 0$ and satisfying the $n-1$ equations in (\ref{eq 2.1}). 
\begin{enumerate}
\item[(i)] Suppose that $D(x)$ is a unit in $RG$. Then $D$ is a $(\sigma, \tau)$-derivation if and only if $\sum_{i=0}^{n-1}x^{iv} = 0$.
\item[(ii)] If $\text{char}(R)$ divides $\text{gcd}(v,n)$, then $D$ is a $(\sigma, \tau)$-derivation. In particular, $\mathcal{D}_{(\sigma, \tau)}(RG)$ is an $R$-module of rank $n$.
\item[(iii)] If $\text{char}(R)$ divides $n$, then $D$ is always a $(\sigma, \sigma)$-derivation. In particular, $\mathcal{D}_{(\sigma, \sigma)}(RG)$ is an $R$-module of rank $n$.
\item[(iv)] Suppose that $n = \text{char}(R)$ is a prime and $D(x)$ is a unit in $RG$. Then $D$ is a $(\sigma, \tau)$-derivation if and only if $\sigma = \tau$.
\item[(v)] Suppose that $\text{char}(R) = 0$, $n$ is a prime, and $D(x)$ is a unit in $RG$. Then $D$ is not a $(\sigma, \tau)$-derivation of $RG$ for any pair $\sigma, \tau$.
\item[(vi)] Suppose that $\text{gcd}(v,n)$ is a unit in $R$. Then $D$ is a $(\sigma, \tau)$-derivation if and only if $\sum_{i=0}^{m-1} c_{i+jd} = 0$ for all $j \in \{0, 1, ..., d-1\}$.
\item[(vii)] Suppose that $n$ is a unit in $R$. Then $RG$ is $(\sigma, \sigma)$-differentially trivial.
\end{enumerate}
\end{corollary}

\begin{corollary}\th\label{corollary 3.3}
Suppose that $\sigma$ and $\tau$ are distinct. Then $RG$ is not $(\sigma, \tau)$-differentially trivial, that is, $\mathcal{D}_{(\sigma, \tau)}(RG) \neq \{0\}$, that is, there always exists a non-zero $(\sigma, \tau)$-derivation of $RG$.
\end{corollary}
\begin{proof}
$v \neq 0$ as $\sigma, \tau$ are distinct. Define $D$ to be an $R$-linear map with $D(1) = 0$, satisfying (\ref{eq 2.1}) and $D(x) = 1-x^{v}$. Then (\ref{eq 3.1}) always holds. 
\end{proof}

We have the following theorem as a consequence of Remark \ref{remark 2.5}, \th\ref{theorem 4.1}, and \th\ref{lemma 3.1}.

\begin{theorem}\th\label{theorem 3.4}
Let $R$ be a commutative unital ring, $G = \langle x \mid x^{n} = 1\rangle$ be a cyclic group of order $n$, and $(\sigma, \tau)$ be a pair of $R$-algebra endomorphisms of $RG$, which are $R$-linear extensions of the group endomorphisms of $G$ (so that $\sigma(x) = x^{u}$ and $\tau(x) = x^{u+v}$ for some $u, v \in \{0, 1, ..., n-1\}$). Then the following statements hold:
\begin{enumerate}
\item[(i)] If $\text{gcd}(v,n)$ is invertible in $R$, then every $(\sigma, \tau)$-derivation of $RG$ is inner, that is, $\mathcal{D}_{(\sigma, \tau)}(RG) = \text{Inn}_{(\sigma, \tau)}(RG)$.
\item[(ii)] If $\text{char}(R)$ divides $\text{gcd}(v,n)$, then $\text{Inn}_{(\sigma, \tau)}(RG) \subsetneq \mathcal{D}_{(\sigma, \tau)}(RG)$, that is, $RG$ has non-zero outer $(\sigma, \tau)$-derivations. 
\end{enumerate}
\end{theorem}

\section{Examples}\label{section 5}
We illustrate \th\ref{theorem 4.1} and \th\ref{theorem 3.4} with some examples.

\begin{example}\th\label{example 4.6}
The derivation defined in the proof of \th\ref{corollary 3.3} is inner with $\beta = x^{-u}$. In fact, given the proof of \th\ref{theorem 4.1}, there exist infinitely many $\beta \in RG$, which make $D$ inner. For example, take $n = 8$ and $(\sigma(x), \tau(x)) = (x^{4}, x^{6})$ so that $v = 2$ and $d=2$. Also, $D(x) = 1-x^{v} = 1-x^{2}$ so that $c_{0} = 1, c_{2} = -1, c_{1} = c_{3} = c_{4} = c_{6} = 0$. Then the augmented matrix $\begin{pmatrix}
A \mid C
\end{pmatrix}$ is row equivalent to the matrix $\left(\begin{array}{cc|c}
I_{n-d} & A'' & C''' \\
0 & 0 & C''
\end{array}\right)$ given by 

$$\left(\begin{array}{cccccccc|c}
1 & 0 & 0 & 0 & 0 & 0 & -1 & 0 & c_{4} \\
0 & 1 & 0 & 0 & 0 & 0 & 0 & -1 & c_{5} \\
0 & 0 & 1 & 0 & 0 & 0 & -1 & 0 & -c_{0}-c_{2} \\
0 & 0 & 0 & 1 & 0 & 0 & 0 & -1 & -c_{1} - c_{3} \\
0 & 0 & 0 & 0 & 1 & 0 & -1 & 0 & -c_{2} \\
0 & 0 & 0 & 0 & 0 & 1 & 0 & -1 & -c_{3} \\
0 & 0 & 0 & 0 & 0 & 0 & 0 & 0 & c_{0} + c_{2} + c_{4} + c_{6} \\
0 & 0 & 0 & 0 & 0 & 0 & 0 & 0 & c_{1} + c_{3} + c_{5} + c_{7}
\end{array}\right).$$

Of course, $\beta = x^{-u} = x^{4}$ makes $D$ inner. However, we have more such $\beta$'s. First, observe that $c_{0} + c_{2} + c_{4} + c_{6} = 0$ and $c_{1} + c_{3} + c_{5} + c_{7} = 0$. Therefore, the system $AB = C$ is consistent with $\beta = \sum_{i=0}^{7} b_{i}x^{i}$, where $b_{0} = b_{6}$, $b_{1} = b_{7}$, $b_{2} = b_{6}$, $b_{3} = b_{7}$, $b_{4} = 1 + b_{6}$, and $b_{5} = b_{7}$, where $b_{6}$ and $b_{7}$ can be chosen freely from $R$.

In fact, observe from above that any $(\sigma, \tau)$-derivation with $D(x) = \sum_{i=0}^{7}c_{i}x^{i}$ ($c_{i} \in R$ for all $0 \leq i \leq 7$) is inner provided $c_{0} + c_{2} + c_{4} + c_{6} = 0$ and $c_{1} + c_{3} + c_{5} + c_{7} = 0$. Further, $\beta \in RG$, which makes $D$ inner (that is, $D(\alpha) = \beta (\sigma - \tau)$ for all $\alpha \in RG$) is given by $\beta = \sum_{i=0}^{7} b_{i}x^{i}$, where $b_{0} = c_{4} + b_{6}$, $b_{1} = c_{5} + b_{7}$, $b_{2} = - c_{0} - c_{2} + b_{6}$, $b_{3} = -c_{1} - c_{3} + b_{7}$, $b_{4} = -c_{2} + b_{6}$, $b_{5} = -c_{3} + b_{7}$, where $b_{6}, b_{7}$ can be chosen freely from $R$.
\end{example}

\begin{example}\th\label{example 4.7}
Suppose that $\text{char}(R) = 2$ and $n=6$. Define $D:RG \rightarrow RG$ to be an $R$-linear map with $D(1) = 0$ and satisfying (\ref{eq 2.1}). Put $D(x) = 1 + x + x^{3} + x^{4}$. Let $(\sigma(x), \tau(x)) = (x^{3}, x^{5})$. Then, using \th\ref{corollary 3.2}, $D$ is a $(\sigma, \tau)$-derivation. Note that $c_{0} = c_{1} = c_{3} = c_{4} = 1$ and $c_{2} = c_{5} = 0$.

$d = \text{gcd}(v,n) = \text{gcd}(2,6) = 2$ and $m = \frac{n}{d} = \frac{6}{2} = 3$. So, in view of \th\ref{theorem 4.1}, $D$ is inner if and only if $c_{0} + c_{2} + c_{4} = 0$ and $c_{1} + c_{3} + c_{5} = 0$. Since these two equations hold here, therefore, $D$ is inner. Now, we determine $\beta$, which makes $D$ inner. The augmented matrix $\begin{pmatrix}
A \mid C
\end{pmatrix}$ is row equivalent to the matrix $\left(\begin{array}{cc|c}
I_{n-d} & A'' & C''' \\
0 & 0 & C''
\end{array}\right)$ given by $$\left(\begin{array}{cccccc|c}
1 & 0 & 0 & 0 & -1 & 0 & c_{3} \\
0 & 1 & 0 & 0 & 0 & -1 & -c_{0} - c_{2} \\
0 & 0 & 1 & 0 & -1 & 0 & -c_{1} \\
0 & 0 & 0 & 1 & 0 & -1 & -c_{2} \\
0 & 0 & 0 & 0 & 0 & 0 & c_{0} + c_{2} + c_{4} \\
0 & 0 & 0 & 0 & 0 & 0 & c_{1} + c_{3} + c_{5} 
\end{array}\right).$$

Therefore, $\beta = \sum_{i=0}^{5} b_{i}x^{i}$, where $b_{0} = 1 + b_{4}$, $b_{1} = 1 + b_{5}$, $b_{2} = 1+b_{4}$, and $b_{3} = b_{5}$; $b_{4}$ and $b_{5}$ can be chosen arbitrarily from $R$. For example, some choices for $\beta$ are: $\beta = 1 + x + x^{2}$, $\beta = x + x^{4}$, $\beta = 1 + x^{2} + x^{3} + x^{5}$, $\beta = x^{3} + x^{4} + x^{5}$.
\end{example}

\begin{example}\th\label{example 3.5}
Suppose $\text{char}(R)$ divides $\text{gcd}(v,n)$. So, according to \th\ref{theorem 3.4}, $RG$ must have non-trivial outer derivations. Define $D:RG \rightarrow RG$ as an $R$-linear map with $D(1) = 0$ and satisfying (\ref{eq 2.1}). Then, by Corollary \ref{corollary 3.2} (ii), $D$ is a $(\sigma, \tau)$-derivation. If we take $D(x)$ to be a trivial unit in $RG$, then by \th\ref{theorem 4.1}, the $(\sigma, \tau)$-derivation $D$ is not inner. Therefore, $D + \text{Inn}_{(\sigma, \tau)}(RG)$ is a non-zero outer derivation of $RG$. Thus, we get a family of non-trivial outer twisted derivations.
\end{example}\vspace{20pt}

\noindent \textbf{\Large Acknowledgements}\vspace{8pt}

\noindent The authors thank and are deeply grateful to the referees and the editor for their critical reviews, comments, and suggestions, which have improved the presentation and quality of the paper. The second author is the ConsenSys Blockchain chair professor. He thanks ConsenSys AG for that privilege.\vspace{30pt}

\noindent \textbf{{\Large Statements and Declarations}}\vspace{8pt}

\noindent \textbf{Ethical Approval}: Not Applicable \vspace{8pt}

\noindent \textbf{Competing Interests}: Not Applicable \vspace{8pt}

\noindent \textbf{Author contributions}: Both authors have contributed equally to the manuscript. \vspace{8pt}

\noindent \textbf{Funding}: Not Applicable \vspace{8pt}

\noindent \textbf{Availability of data and materials}: Not Applicable

\bibliographystyle{plain}
\end{document}